\documentclass[12pt, amstex]{article}

\usepackage{graphicx}

\usepackage{amsfonts}
\usepackage{float}
\usepackage{flafter}

\usepackage{color}

\labelwidth\leftmargini\advance\labelwidth-\labelsep

\def\R{\relax\ifmmode I\!\!R\else$I\!\!R$\fi}

\def\Z{\relax\ifmmode Z\!\!\!Z\else$Z\!\!\!Z$\fi}

\def\C{\relax\ifmmode C\!\!\!\!I\else$C\!\!\!\!I$\fi}

\def\K{\relax\ifmmode I\!\!K\else$I\!\!K$\fi}

\def\N{\relax\ifmmode I\!\!N\else$I\!\!N$\fi}

\newcounter{defcounter}[section]
{\vspace{0.1cm}\begin{sloppypar}\noindent\stepcounter{defcounter}{\bfseries
Definition
      \thesection.\thedefcounter}}%
{\end{sloppypar}\vspace{0.1cm}}

\newtheorem{corollary}{Corollary}[section]
\newtheorem{lemma}{Lemma}[section]
\newtheorem{theorem}{Theorem}[section]

\newtheorem{proposition}{Proposition}[section]

\newcommand{\proof}{{\bf Proof.} }

\newcommand{\qed}{\hfill $\square$}

\begin{document}
\thispagestyle{empty}
\begin{center}
{\Large {\bf Nonuniform Markov geometric measures}}
\end{center}
\begin{center}
J. Neunh\"auserer\\
University of Applied Science Berlin \\
neunchen@aol.com
\end{center}
\begin{center}
\begin{abstract}
We generalize results of Fan and Zhang \cite{[FZ]} on absolute continuity and singularity of the golden Markov geometric series to nonuniform stochastic series given by arbitrary Markov process. In addition we describe an application of these results in fractal geometry. 
\end{abstract}
\end{center}
\section{Introduction}
Geometric probability measures on the real line induced by Bernoulli processes are intensely studied in geometric measure theory since the work of Erd\"os \cite{[ER1],[ER2]}. We have a couple of results on singularity and dimension resp. absolute continuity and density of these measures. We refer here to Peres, Schlag and Solomyak \cite{[PSS]} for a nice overview on results about Bernoulli convolutions and to \cite{[NE1],[NE2],[NW]} for results on nonuniform self self-similar measures given Bernoulli processes.~\\
 To achieve further progress, it seems natural to drop the assumption of independence of the process and to study the properties Markov geometric measures. A result in this direction is due to Fan and Zhang \cite{[FZ]}. They proved that the distribution of random power series
\[ \sum_{i=1}^{\infty}X_{i}\beta^{i}\]
for the golden Markov process $(X_{i})$ is singular, if $\beta\in(0,\frac{\sqrt{5}-1}{2})$ and absolutely continuous for almost all $\beta\in (\frac{\sqrt{5}-1}{2},0.739)$. In this paper we generalize this approach to nonuniform stochastic series
\[ \sum_{i=1}^{\infty}X_{i}\beta_{1}^{\hat X_{i}}\beta_{0}^{i-\hat X_{i}}\]
for arbitrary 2-step Markov process $(X_{i})$. We find an upper bound on the Hausdorff dimension of the measures which gives us a domain of singularity, see Theorem 3.1. Outside this domain we prove generic absolute continuity of the measure under the condition that the parameter values $\beta_{0}$, $\beta_{1}$ are below $0.739$. In addition we find a subdomain where the measures generically have a density in $L^2$, see Theorem 3.2.  The bound $0.739$ that appears here is due to the transversality techniques used in the proofs, see Section 5. Our result on absolute continuity can be applied in fractal geometry to prove a weighted version of the classical Moran formula \cite{[MO]} for the Hausdorff dimension of certain self affine sets, see Theorem 3.3. \\
The rest of the paper is organized as follows: In the next section we introduce three different descriptions of Markov geometric measures and in Section 3 we formalize our result recalling some notions in geometric measure theory. In Section 4 the reader find the proof of the result on singularity, the next section contains the proof of results on absolute continuity and in the last section we prove our result in fractal geometry.
\section{Description of the measures}
For $p\in(0,1)$ let $(X_{i})_{i \in \mathbb{N}}$ be the Markov chain with two states $0$ and $1$, given by the transition matrix
\[ M=\left(\begin{array}{cc} p_{00}  & p_{01} \\p_{10}&p_{11} \end{array}\right)=\left(\begin{array}{cc} p  &1-p \\1&0 \end{array}\right)\]
and initial probability
\[ (p_{0},p_{1})=(\frac{1}{2-p},\frac{1-p}{2-p}).\]
Consider random power series in two variables $\beta_{0},\beta_{1}\in(0,1)$, given by
\[ S=\sum_{i=1}^{\infty}X_{i}\beta_{1}^{\hat X_{i}}\beta_{0}^{i-\hat X_{i}}, \]
where the random variable $\hat X_{i}$ counts the number of entries in the chain that are one,
\[ \hat X_{i}= \sum_{k=1}^{i}X_{k}.\]
We call the distribution of this random power series,
\[ \mu(B)=\mu_{\beta_{0},\beta_{1},p}(B)=Prob(S\in B)\]
for Borel sets $B\subseteq \mathbb{R}$, a nonuniform Markov geometric measure on the real line. If $\beta_{0}=\beta_{1}$ and $p$ is the inverse of the golden mean, the measures defined here are exactly the distribution of
golden Markov geometric series introduced by Fan and Zhang \cite{[FZ]}. To get results on the properties of the measures we introduce another descriptions.~\\
Given the matrix
\[ A=\left(\begin{array}{cc} a_{00}  & a_{01} \\a_{10}&a_{11} \end{array}\right)=\left(\begin{array}{cc} 1  &1 \\1&0 \end{array}\right),\]
consider the Fibonacci subshift
\[\Sigma_{A}=\{(s_{k})| a_{s_{i}s_{i+1}}=1\}\subseteq\{0,1\}^{\mathbb{N}}\]
with the natural product metric
\[ d((s_{k}),(t_{k}))=\sum_{k=1}^{\infty} |s_{k}-t_{k}|2^{-k}.\]
The Markov chain  $(X_{i})_{i \in \mathbb{N}}$ induces a Borel probability measure $\nu_{p}$ supported by $\Sigma_{A}$. We introduce a coding map $\pi:\Sigma_{A}\longmapsto\mathbb{R}$ by
\[ \pi((s_{k}))=\sum_{i=1}^{\infty}s_{i}\beta_{1}^{\sharp_{i}((s_{k})) }\beta_{0}^{i-\sharp_{i}((s_{k}))},\]
where $\sharp_{i}((s_{k}))=\sum_{k=1}^{i}s_{k}$. Obviously the Markov geometric measure $\mu$ is given by
\[ \mu(B)=\pi(\nu_{p})(B)=\nu_{p}(\pi^{-1}(B))\]
for Borel sets $B\subseteq \mathbb{R}$. Note that $\mu$ is supported by $\pi(\Sigma_{A})$, which is a Cantor set in the case $\beta_{0}+\beta_{0}\beta_{1}<1$ and contains an interval if  $\beta_{0}+\beta_{0}\beta_{1}\ge 1$. Furthermore note that the set $(\pi_{\beta_{0},\beta_{1}}\times\pi_{\tau_{0},\tau_{1}})(\Sigma_{A})\subseteq \mathbb{R}^2$
for $\beta_{0},\beta_{1},\tau_{0},\tau_{1}\in(0,1)$ is a self affine fractal, modeled by the subshift $\Sigma_{A}$, see chapter 5 of \cite{[PE]}. It supports the Borel probability measure $(\pi_{\beta_{0},\beta_{1}}\times\pi_{\tau_{0},\tau_{1}})(\nu)$
on $\mathbb{R}^2$.
\section{Results}
For convenience recall some notions from measure theory, see also \cite{[BI]} or \cite{[MA]}. A Borel probability measure $\mu\in\mathfrak{M}$  is totally singular with respect to the Lebesgue measure $\mathfrak{L}$, if there is a Borel set $B$ of Lebesgue measure zero, $\mathfrak{L}(B)=0$, with
$\mu(B)=1$. The measure is absolutely continuous with respect to the Lebesgue measure, if $\mu(B)>0$ holds only for Borel sets $B$ of positive Lebesgue measure, $\mathfrak{L}(B)>0$. Moreover any measure $\mu$ my be decomposed into singular $\mu_{S}$ and absolute continuous part $\mu_{A}$, which means $\mu=\mu_{S}+\mu_{A}$. It follows from \cite{[JW]} that the distributions of random power series are of pure type, either absolutely continuous ($\mu_{S}=0$) or totally singular ($\mu_{A}=0$) with respect to the Lebesgue measure.~\\
In order to state our result on singularity of Markov geometric measure we define the Hausdorff dimension of a Borel probability measure $\mu$ on the line by
\[ \dim_{H}\mu=\inf\{\dim_{H}A|\mu(A)=1\},\]
where $\dim_{H}A$ is the Hausdorff dimension of a subset $A$ of $\mathbb{R}$. We refer to the book of Falconer \cite{[FA]} or the book of Pesin \cite{[PE]} for an introduction to dimension theory. Note that $\dim_{H}\mu<1$ obviously implies the singularity of $\mu$ with respect to the Lebesgue measure. We will prove the following theorem in Section 4:
\begin{theorem} For all $\beta_{0},\beta_{1},p\in(0,1)$ with
\[ \beta_{0}\beta_{1}^{1-p}<p^{p}(1-p)^{1-p}\]
the Markov geometric measure $\mu_{\beta_{0},\beta_{1},p}$ is singular with
\[ \dim_{H}\mu_{\beta_{0},\beta_{1},p}\le \frac{p\log(p)+(1-p)\log(1-p)}{\log(\beta_{0})+(1-p)\log(\beta_{1})}.\]
\end{theorem}~\\
We note that the equality in the dimension estimate above follows from classical results in the case that $\mu$ is supported by a Cantor set, that is  means $\beta_{0}+\beta_{0}\beta_{1}<1$, see \cite{[PE]}.~\\~\\
To state results on absolute continuity, recall that by the theorem of Radon-Nikodym a Borel probability measure $\mu$ is absolutely continuous, if and only if it has a density $f\in L^{1}$
\[ \mu=\int f d\mathfrak{L}.\]
For an absolutely continuous measure we may thus ask if this density is $L^2$, which means
\[ \int f^2 d \mathfrak{L}<\infty.\]
In Section 5 we will prove the following result:
\begin{theorem}
Fix $p\in(0,1)$. For almost all  $\beta_{0},\beta_{1}\in(0,0.739)$ with
 \[ \beta_{0}\beta_{1}^{1-p}>p^{p}(1-p)^{1-p}\]
the Markov geometric measure $\mu_{\beta_{0},\beta_{1},p}$ is absolutely continuous and has density in $L^2$ for almost all $\beta_{0},\beta_{1}\in(0,0.739)$ with
\[ (\beta_{0}-p^2)\beta_{1}>(1-p)^2.\]
\end{theorem}
Note that the lower bound on absolute continuity here is sharp by Theorem 3.2. The upper bound is due to the techniques we will use in the proof and is not expected to be be sharp. ~\\In the last section we will apply our result on absolute continuity of Markov geometric measure to obtain the following result in fractal geometry:
\begin{theorem}
For all $\tau_{1},\tau_{2}\in(0,1)$ with $\tau_{1}+\tau_{1}\tau_{2}<1$ and almost all $\beta_{1},\beta_{2}\in(0,0.739)$ with $\beta_{1}+\beta_{1}\beta_{2}>1$, the Hausdorff dimension $d$ of the self affine fractal
\[\Lambda=(\pi_{\beta_{1},\beta_{2}}\times\pi_{\tau_{1},\tau_{2}})(\Sigma_{A})\]
is given by the unique solution $d>1$ of
\[ \beta_{1}\tau_{1}^{d-1}+\beta_{1}\beta_{2}(\tau_{1}\tau_{2})^{d-1}=1.\]
\end{theorem}~\\
We remark that this is a weighted version of the classical Moran formula $\tau_{1}^{d}+(\tau_{1}\tau_{2})^{d}=1$
for the Hausdorff dimension of the self-similar sets $\pi_{\tau_{1},\tau_{2}}(\Sigma_{A})$ modeled by the subshift, see  \cite{[MO]} and  \cite{[PE]}.
\section{Singularity}
Given $\beta_{0},\beta_{1}\in(0,1)$, we define an adapted metric on the Fibonacci subshift $\Sigma_{A}$ by
\[ d_{\beta_{0},\beta_{1}}((s_{k}),(t_{k}))= \beta_{1}^{\sharp_{i}((s_{k}))}\beta_{0}^{i-\sharp_{i}((s_{k})) },\]
where
\[ i=i((s_{k}),(t_{k}))=\max\{k|s_{k}=t_{k}\} \]
and
\[ \sharp_{i}((s_{k}))=\sum_{k=1}^{i}s_{k}.\]
From \cite{[NE1]} we know that the following lemma holds.
\begin{lemma}
The coding map $\pi:\Sigma_{A}\longmapsto\mathbb{R}$ is Lipschitz continuous with respect to the adapted metric $d_{\beta_{0},\beta_{1}}$.
\end{lemma}
Now we need  the following proposition.
\begin{proposition}
For $p\in(0,1)$ let $\nu$ be the Borel probability measure on $\Sigma_{A}$, given by the matrix
\[\left(\begin{array}{cc} p  &1-p \\1&0 \end{array}\right)\]
and the  vector
\[ (\frac{1}{2-p},\frac{1-p}{2-p}).\]
The Hausdorff dimension of $\nu$  with respect to the adapted metric $d_{\beta_{0},\beta_{1}}$ on $\Sigma_{A}$ is given by
\[  \dim_{H}\nu_{p}= \frac{p\log(p)+(1-p)\log(1-p)}{\log(\beta_{0})+(1-p)\log(\beta_{1})}.\]
\end{proposition}
\proof
We know from Parry \cite{[PA]} that the entropy of the measure $\nu_{p}$ is given by
\[h(\nu_{p})=-(\frac{p}{2-p}\log(p)+\frac{1-p}{2-p}\log(1-p)).\]
Moreover by Shannon's local entropy theorem for $\nu_{p}$-almost all sequences $(s_{k})$ in $\Sigma_{A}$
\[\lim_{n\longmapsto\infty}-\frac{1}{n}\log\nu_{p}([s_{1},s_{2},\dots, s_{n}])=h(\nu_{p}),\]
where $[s_{1},s_{2},\dots s_{n}]$ denotes the cylinder set. The diameter of this cylinder set with respect to the metric $d_{\beta_{0},\beta_{1}}$ is given by
\[ |[s_{1},s_{2},\dots ,s_{n}]|=\beta_{1}^{\sharp_{n}((s_{k}))}\beta_{0}^{n-\sharp_{n}((s_{k})) },\]
and by Birkhoff's ergodic theorem we have for $\nu_{p}$-almost all sequences $(s_{k})$ in $\Sigma_{A}$
\[ \lim_{n\longmapsto\infty}-\frac{1}{n}\log(\beta_{1}^{\sharp_{n}((s_{k}))}\beta_{0}^{n-\sharp_{n}((s_{k}))})=-(\frac{p}{2-p}\log(\beta_{0})+\frac{1-p}{2-p}\log(\beta_{1}))\]
for almost all sequences $(s_{k})$ in $\Sigma_{A}$. Hence the local dimension of the measure $\nu$ is given by
\[ d(\nu_{p},(s_{k}))=\lim_{n\longmapsto\infty}\frac{\log(\nu_{p}([s_{1},s_{2},\dots s_{n}]))}{ \log(|[s_{1},s_{2},\dots s_{n}]|)}= \frac{p\log(p)+(1-p)\log(1-p)}{\log(\beta_{0})+(1-p)\log(\beta_{1})}\]
almost surely and our result follows from the local mass distribution principle, see \cite{[FA]}.\qed~~\\~\\
We are now prepared to prove theorem 3.2.\\~\\
{\bf Proof of Theorem 3.2} It is well known \cite{[FA]}, that Lipschitz maps do not increase Hausdorff dimension.  Hence $\dim_{H}\mu_{\beta_{0},\beta_{1},p}\le\dim_{H}\nu$ since $\mu_{\beta_{0},\beta_{1},p}=\pi(\nu)$ with respect to the metric  $d_{\beta_{0},\beta_{1}}$ by lemma 3.1. Now proposition 3.2 immediately gives the dimension estimate in theorem 3.2. For the singularity assertion in theorem 3.2 just note that the dimension estimate
implies $\dim_{H}\mu_{\beta_{0},\beta_{1},p}<1$ if
$ \beta_{0}\beta_{1}^{1-p}<p^{p}(1-p)^{1-p}$. \qed
\section{Absolute continuity}
In order to prove results on absolute continuity of the Markov geometric measures we use transversality techniques developed by Peres and Solomyak \cite{[PS1],[PS2]} based on an idea by Pollicott and Simon \cite{[PS]}. For a set $\mathfrak{B}$ of real analytic functions we say that the $\delta$-transversality condition holds on an interval $I$, if
\[ f(x)<\delta \Rightarrow  f'(x)<-\delta\]
holds for all $x\in I$ and $f\in\mathfrak{B}$ where $\delta>0$ is independent of $x$ and $f$. This means that the graph of the function
$f$ crosses all horizontal lines that it meets below height $\delta$  transversely with slope
at most $-\delta$. In \cite{[PS1]} it is proved that there is an $\delta$ such that the  $\delta$-transversality condition holds for
 \[ \mathfrak{B}=\{f|f(x)=1+\sum_{k=1}a_{k}x^{k},\quad |a_{k}|\le 1\}\]
on the interval $[b,0.649]$ for any $b>0$. Fan and Zhang \cite{[FZ]} succeeded in extending the upper bound here a little bit if the coefficients of the power series come from a Markov shift $\Sigma_{A}$. We state the result as a lemma.
\begin{lemma}
For any $b>0$ there is a $\delta$ such that the  $\delta$-transversality condition holds for the set of power series
 \[ \mathfrak{B}_{A}=\{f|f(x)=1+\sum_{k=1}^{\infty}(a_{k}s_{k}-b_{k}t_{k})x^{k},\quad a_{k},b_{k}\in(0,1],\quad(s_{k}),(t_{k})\in\Sigma_{A}\}\]
on the interval $[b,0.739]$.
\end{lemma}
In fact this lemma is proved in \cite{[FZ]} for $a_{k}=b_{k}=1$. But as remarked by Peres and Solomyak, see lemma 5.1 of \cite{[PS2]}, the technique using (*)-functions generalizes to the case of smaller coefficients and this technique is applied in \cite{[FZ]}.  Using transversality lemma we will prove the following lemma, which will be essential.
\begin{lemma}For all $\underline{\beta}>0$ and $c\in(0,1)$ there is constant $D>0$ such that for all sequences $(s_{k}),(t_{k})\in\Sigma_{A}$
\[\mathfrak{L}(\{
\beta\in[\underline{\beta},0.739]     |~|\pi_{\beta,c\beta}((s_{k}))-\pi_{\beta,c\beta}((t_{k}))|\le r   \})\le Dr\underline{\beta}^{-i}c^{-\sharp_{i}((s_{k}))}\]
where $i=i((s_{k}),(t_{k}))=\max\{k|s_{k}=t_{k}\}$.
\end{lemma}
\proof Define $i$ depending on the sequences $(s_{k}),(t_{k})\in\Sigma_{A}$ as above and suppress these indices in the following. Using the definition of the coding map $\pi$ in section 2 we obtain,
\[ \phi_{(s_{k}),(t_{k})}(\beta)=\pi_{\beta,c\beta}((s_{k}))-\pi_{\beta,c\beta}((t_{k}))=\sum_{n=1}^{\infty}\left(s_{n}c^{\sharp_{n}((s_{k}))}-t_{n}c^{\sharp_{n}((t_{k}))}\right)\beta^{n}\]
\[= \beta^{i}\sum_{n=1}^{\infty}\left(s_{i+n}c^{\sharp_{i+n}((s_{k}))}-t_{i+n}c^{\sharp_{i+n}((t_{k}))}\right)\beta^{n}\]
\[ =\beta^{i+1}\left(s_{i+1}c^{\sharp_{i+1}((s_{k}))}-t_{i+1}c^{\sharp_{i+1}((t_{k}))}\right)(1+\sum_{n=1}^{\infty}\frac{\left(s_{i+n+1}c^{\sharp_{i+n+1}((s_{k}))}-t_{i+n+1}c^{\sharp_{i+n+1}((t_{k}))}\right)}{\left(s_{i+1}c^{\sharp_{i+1}((s_{k}))}-t_{i+1}c^{\sharp_{i+1}((t_{k}))}\right) }\beta^{n})\]
\[ =\beta^{i+1}c^{\sharp_{i}((s_{k}))+1}(\bar{s}_{1}-\bar{t}_{1})\underbrace{\left(1+\sum_{n=1}^{\infty}\frac{\bar{s}_{n+1}c^{\sharp_{n+1}((\bar{s}_{k}))-1}-\bar{t}_{n+1}c^{\sharp_{n+1}((\bar{t}_{k}))-1}}{\bar{s}_{1}-\bar{t}_{1}}\beta^{n}
\right)}_{:=\psi_{(s_{k}),(t_{k})}(\beta)},\]
where $(\bar{s}_{k})$ is the $i$-times shifted sequence $(s_{i+1},s_{i+2},\dots)\in\Sigma_{A}$ and accordingly $(\bar{t}_{k})$ is the sequence $(t_{i+1},t_{i+2},\dots)\in\Sigma_{A}$.
Note that if $\bar{s}_{n+1}=1$ we have $c^{\sharp_{n+1}((\bar{s}_{k}))-1}\in(0,1)$, the same holds for $\bar t$. Hence we see (exchanging the role of $t$ and $s$ if $\bar{s}_{1}-\bar{t}_{1}=-1$) that the power series $\psi_{(s_{k}),(t_{k})}$ here fall into the class $\mathfrak{B}_{A}$, defined in lemma 5.1. Applying $\rho$-transversality we get
\[\mathfrak{L}\{\beta\in[\underline{\beta},0.739]   ~  |~|\phi_{(s_{k}),(t_{k})}(\beta)|\le r   \}\]
\[=\mathfrak{L}\{\beta\in[\underline{\beta},0.739]   ~  |~|\psi_{(s_{k}),(t_{k})}(\beta)|\le r\beta^{-(i+1)}c^{-(\sharp_{i}((s_{k}))+1)}\}  \le 2\rho^{-1}\beta^{-(i+1)}c^{-(\sharp_{i}((s_{k}))+1)}r.\]
Thus our lemma holds with $D=2\rho^{-1}\underline{\beta}^{-1}c^{-1}$.
\qed~\newline\newline
Now we are ready to prove our main result:
\begin{theorem}
Let $p\in(0,1)$, $q\in(1,2]$ and $c\in(0,1]$. The density of the measures $\mu_{\beta,c\beta,p}$ is in $L^{q}$
for almost all $\beta\in(0,0.739)$ with
\[ \left(\frac{p^{q}}{\beta^{q-1}}+\frac{(1-p)^{q}}{(c\beta^2)^{q-1}}\right)<1.\]
\end{theorem}
\proof Fix $p$, $q$ and $c$ during the proof.
Define the (lower) local density of a measure $\mu$ on the real line
by
\[ \underline{D}(\mu,x)=\liminf_{r\longrightarrow 0}
\frac{\mu(B_{r}(x))}{2r}.\]
If
\[ \int(\underline{D}(\mu,x))^{q-1}d\mu(x)<\infty,\]
then $\mu$ is absolute continuous and has
density in $L^{q}$ by \\S 2.12 of \cite{[MA]}.
Thus it suffices to show that
\[\Im(\underline{\beta}):=\int_{\underline{\beta}}^{0.739}\int~ (D(\mu_{p,\beta,c\beta},x))^{q-1}d\mu_{p,\beta,c\beta}(x)~d\beta<\infty \]
holds for all $\underline{\beta}>\underline{\beta}(p,q,c)$.
Applying Fatou's lemma and changing the order of integration and variables using $\mu_{p,\beta,c\beta}=\pi_{\beta,c\beta}(\nu_{p})$
and estimating with H\"older inequality  $\int f^{q-1}\le C(\int f)^{q-1}$, we obtain
\[ \Im(\underline{\beta})\le    \liminf_{r\longrightarrow 0}\frac{1}{(2r)^{q-1}}
\int_{0}^{0.739}\int_{\Sigma_{A}}~\left(\mu_{p,\beta,c\beta}(B_{r}(x))\right)^{q-1}d\mu_{p,\beta,c\beta}dx\]
\[ = \liminf_{r\longrightarrow 0}\frac{1}{(2r)^{q-1}}
\int_{\Sigma_{A}}\int_{0}^{0.739}~\left(\mu_{p,\beta,c\beta}(B_{r}(\pi_{\beta,c\beta}((s_{k})))\right)^{q-1}dx~d\nu_{p}\]
\[ \le C\liminf_{r\longrightarrow 0}\frac{1}{(2r)^{q-1}}
\int_{\Sigma_{A}}\left(\int_{0}^{0.739}~\mu_{p,\beta,c\beta}(B_{r}(\pi_{\beta,c\beta}((s_{k})))dx\right)^{q-1}~d\nu_{p}\]
\[= C~\liminf_{r\longrightarrow 0}\frac{1}{(2r)^{q-1}}
\int_{\Sigma_{A}}\left(\int_{\Sigma_{A}}\mathfrak{L}(\{
\beta\in[\underline{\beta},0.739]     |~|\pi_{\beta,c\beta}((s_{k}))-\pi_{\beta,c\beta}((t_{k}))|\le r   \})~d\nu_{p}\right)^{q-1}d\nu_{p}
.\]
Using lemma 5.2 we continue with
\[\Im(\underline{\beta})\le  C~\liminf_{r\longrightarrow 0}\frac{1}{(2r)^{q-1}}\int_{\Sigma_{A}}\left(\int_{\Sigma_{A}}~Dr\underline{\beta}^{-i}c^{-\sharp_{i}((s_{k}))}d\nu_{p}\right)^{q-1}d\nu_{p}\]
\[ =C(D/2)^{q-1}\int_{\Sigma_{A}}\left(\int_{\Sigma_{A}}\underline{\beta}^{-i}c^{-\sharp_{i}((s_{k}))}d\nu_{p}\right)^{q-1}d\nu_{p}\]
\[ \le C(D/2)^{q-1}\max\{p,1-p\}\int_{\Sigma_{A}}\left(\sum_{n=1}^{\infty}\underline{\beta}^{-n}c^{-\sharp_{n}((s_{k}))}\nu_{p}([s_{1},\dots,s_{n}])\right)^{q-1}d\nu_{p}\]
\[ \le C(D/2)^{q-1}\max\{p,1-p\}\sum_{n=1}^{\infty}\int_{\Sigma_{A}}\underline{\beta}^{-n(q-1)}c^{-\sharp_{n}((s_{k}))(q-1)}\nu_{p}([s_{1},\dots,s_{n}])^{q-1}d\nu_{p}\]
\[ =  \ C(D/2)^{q-1}\max\{p,1-p\}\sum_{n=1}^{\infty}\sum_{[s_{1},\dots,s_{n}]\subseteq\Sigma_{A}}\underline{\beta}^{-n(q-1)}c^{-\sharp_{n}((s_{k}))(q-1)} \nu_{p}([v_{1},\dots,v_{n}])^{q}\]
where the sum is taken over cylinder sets in $\Sigma_{A}$. We now write the summation here in another form using the structure of $\Sigma_{A}$. For a sequence $v=(v_{1},\dots,v_{n})\in\{1,(-1,1)\}^{n}$ let $\sharp_{1}(v)$ be the number of entries in $v$ that are one and by  $\sharp_{-1,1}(v)$ be the number of entries that are $(-1,1)$ one. With this notation the measure of the corresponding cylinder set is given by
\[ \nu_{p}([v_{1},\dots,v_{n}])= \nu_{p}([v_{1}])p^{\sharp_{1}(v)}(1-p)^{\sharp_{-1,1}(v)}\le \max\{\frac{1}{2-p},\frac{1-p}{2-p}\}p^{\sharp_{1}(v)}(1-p)^{\sharp_{-1,1}(v)}\]
and the length of the cylinder is obviously $\sharp_{1}(v)+2\sharp_{-1,1}(v)$. Hence $\Im(\underline{\beta})$ is up to the constant
\[ C(D/2)^{q-1}\max\{p,1-p\}\max\{\frac{1}{2-p},\frac{1-p}{2-p}\}\]
bounded from above by
\[ \sum_{n=1}^{\infty}\sum_{v\in\{1,(-1,1)\}^{n}}\underline{\beta}^{-(\sharp_{1}(v)+2\sharp_{-1,1}(v))(q-1)}c^{-(q-1)\sharp_{-1,1}(v)} p^{q(\sharp_{1}(v)}(1-p)^{q\sharp_{-1,1}(v)}\]
\[ = \sum_{n=1}^{\infty}\sum_{v\in\{1,(-1,1)\}^{n}}(\underline{\beta}^{1-q}p^{q})^{-\sharp_{1}(v)}(c^{1-q}\underline{\beta}^{2(1-q)}(1-p)^{q})^{\sharp_{-1,1}(v)}\]
\[ =  \sum_{n=1}^{\infty} \left(\underline{\beta}^{1-q}p^{q}+(c\underline{\beta}^{2})^{1-q}(1-p)^{q}\right)^{n}.\]
The geometric series here converges if  $\underline{\beta}^{1-q}p^{q}+(c\underline{\beta}^{2})^{1-q}(1-p)^{q}<1$, which is our assumption.  \qed~\\~\\
Our theorem has the following corollary on absolute continuity.
\begin{corollary}
For all $p\in(0,1)$ and $c\in(0,1]$ the measure $\mu_{\beta,c\beta,p}$ is absolutely continuous
for almost all $\beta\in(0,0.739)$ with
\[ (\frac{p}{\beta})^{p}(\frac{1-p}{c\beta^{2}})^{1-p}< 1.\]
\end{corollary}
\proof
By Theorem 5.1 the measure $\mu_{\beta,c\beta,p}$ is absolutely continuous for all $q\in (0,1]$ and almost all
$\beta\in(0,0.739)$ with
\[ \left(p^{q}+\frac{(1-p)^{q}}{(c\beta)^{q-1}}\right)^{1/(q-1)}<\beta.\]
Taking the limit $q\longmapsto 1$ and using the rule of L'Hospital we obtain that $\mu_{\beta,c\beta,p}$ is absolutely continuous
for almost all $\beta\in(0,0.739)$ with
\[     p^{p}(\frac{1-p}{c\beta})^{1-p}  < \beta,\]
which is equivalent to the condition stated in the corollary.  \qed~\\~\\
{\bf Proof of theorem 3.3} By the the theorem of Fubini the last corollary is stronger than our assertion on absolute continuity in theorem 3.3. Setting $q=2$ in theorem 5.1 and using Fubini again, we get the assertion on $L^2$ density in theorem 3.3. \qed
\section{An application to fractal geometry}
Our aim here is to prove Theorem 3.4. on the Hausdorff dimension of the sets
\[\Lambda(\beta_{0},\beta_{1},\tau_{0},\tau_{1})=(\pi_{\beta_{0},\beta_{1}}\times\pi_{\tau_{0},\tau_{1}})(\Sigma_{A}),\]
introduced in Section two. From the view of fractal geometry the attractor of the iterated function system $G_{0},G_{1}:\mathbb{R}^2\rightarrow \mathbb{R}^2$ with
\[ G_{0}(x,y)=(\beta_{0}x,\tau_{0}y) \quad  G_{1}(x,y)=(\beta_{1}x+\beta_{1},\tau_{1}y+\tau_{1})\]
is given by $(\pi_{\beta_{0},\beta_{1}}\times\pi_{\tau_{0},\tau_{1}})(\{0,1\}^{\mathbb{N}})$ and $\Lambda$ is the subset of this attractor
modeled by $\Sigma_{A}$, see chapter 5 of \cite{[PE]}. By a standard covering argument from fractal geometry (see \cite{[PW]}) we have get following proposition:
\begin{proposition}
For all $\tau_{0},\tau_{1},\beta_{0},\beta_{1}\in(0,1)$ with $\tau_{0}+\tau_{0}\tau_{1}<1$ and $\beta_{0}+\beta_{0}\beta_{1}\ge 1$ we have
\[ \dim_{H}\Lambda(\beta_{0},\beta_{1},\tau_{0},\tau_{1})\le d ,\]
where $d$ is given by the unique solution $d>1$ of
\[ \beta_{0}\tau_{1}^{d-1}+\beta_{0}\beta_{1}(\tau_{0}\tau_{1})^{d-1}=1.\]
\end{proposition}
In fact the box-counting dimension of the set is given by $d$ and the box-counting dimension is an upper bound on the Hausdorff dimension in general. To show that $d$ generically is a lower bound, we use the measures
\[ \hat\mu(\beta_{0},\beta_{1},\tau_{0},\tau_{1},p):=(\pi_{\beta_{0},\beta_{1}}\times\pi_{\tau_{0},\tau_{1}})(\nu_{p}),\]
that are supported on $\Lambda(\beta_{0},\beta_{1},\tau_{0},\tau_{1})$. We can express the Hausdorff dimension of these measures in terms of the Hausdorff dimension of the geometric Markov measures $\mu_{\beta_{0},\beta_{1}}=\pi_{\beta_{0},\beta_{1}}(\nu)$.
\begin{proposition}
For all $\tau_{0},\tau_{1},\beta_{0},\beta_{1}\in(0,1)$ with $\tau_{0}+\tau_{0}\tau_{1}<1$ and $\beta_{0}+\beta_{0}\beta_{1}\ge 1$, we have
\[ \dim_{H}\hat\mu(\beta_{0},\beta_{1},\tau_{0},\tau_{1},p)=\frac{p\log(p)+(1-p)\log(1-p)}{\log(\tau_{0})+(1-p)\log(\tau_{1})}\] \[+(1-\frac{\log(\beta_{0})+(1-p)\log(\beta_{1})}{\log(\tau_{0})+(1-p)\log(\tau_{1})})\dim_{H}\mu_{\beta_{0},\beta_{1},p}.\]
\end{proposition}
This proposition can be proved using the dimension theory of dynamical systems, see chapter 8 of \cite{[PE]}, exactly in the same way as theorem 2 of \cite{[NE1]}. Just replace the Bernoulli measures there by Markov measures and use the formula of Parry for the entropy of these measures, compare section 4. The proposition may be also proved using the general dimension theory of iterated function systems \cite{[FH]}, by relating the conditional entropies that appear in this work to the dimension of projected measures. Using both proposition above it is easy the prove our last result.~\\~\\
{\bf Proof of theorem 3.4} Under the assumption Proposition 6.1 and Proposition 6.2 let
$d>1$ be the solution of
$ \beta_{0}\tau_{0}^{d-1}+\beta_{0}\beta_{1}(\tau_{0}\tau_{1})^{d-1}=1$ and let $p=\beta_{0}\tau_{0}^{d-1}$. By Proposition 6.2 we have
$\dim_{H}\hat\mu(\beta_{0},\beta_{1},\tau_{0},\tau_{1},p)=d$ and by Proposition 6.1 we get $\dim_{H}\Lambda(\beta_{0},\beta_{1},\tau_{0},\tau_{1})=d$, if $\dim_{H}\mu_{\beta_{0},\beta_{1},p}=1$. It remains to show that
for all $\tau_{0},\tau_{1}\in(0,1)$ with $\tau_{0}+\tau_{0}\tau_{1}<1$ and almost all $\beta_{0},\beta_{1}\in(0,0.739)$ the measure $\mu_{\beta_{0},\beta_{1},p}$ for $p=\beta_{0}\tau_{0}^{d-1}$ is absolutely continuous and hence has dimension one. Note that
\[ \frac{p^2}{\beta_{0}}+ \frac{(1-p)^2}{\beta_{0}\beta_{1}}=\beta_{0}\tau_{0}^{2(d-1)}+\beta_{0}\beta_{1}(\tau_{0}\tau_{1})^{2(d-1)}<1,\]
hence the statement follows from theorem 5.1 for $q=2$.  \qed

\end{document}